\documentclass[11pt,reqno]{amsart}
\usepackage{a4wide}
\usepackage{euscript,amsmath,amssymb,amsbsy,mathabx}
\usepackage[arrow,matrix,curve]{xy}
\usepackage{array,longtable,graphicx,enumitem,url}

\usepackage{xcolor}
\usepackage[colorlinks=true,linkcolor=black,urlcolor=black,citecolor=black]{hyperref}

\numberwithin{equation}{section}\numberwithin{figure}{section}

\newcounter{msct}[section]\renewcommand{\themsct}{\thesection.\arabic{msct}}
\newenvironment{m-theorem}{\vskip3pt\refstepcounter{msct}\trivlist \itemindent 0pt %
\item[\hskip\labelsep\bf Theorem \themsct]\it\ignorespaces}{\endtrivlist\vskip2pt}
\newenvironment{m-proposition}{\vskip3pt\refstepcounter{msct}\trivlist \itemindent 0pt %
\item[\hskip\labelsep\bf Proposition \themsct]\it\ignorespaces}{\endtrivlist\vskip2pt}
\newenvironment{m-corollary}{\vskip3pt\refstepcounter{msct}\trivlist \itemindent 0pt %
\item[\hskip\labelsep\bf Corollary \themsct]\it\ignorespaces}{\endtrivlist\vskip2pt}
\newenvironment{m-lemma}{\vskip3pt\refstepcounter{msct}\trivlist \itemindent 0pt %
\item[\hskip\labelsep\bf Lemma \themsct]\it\ignorespaces}{\endtrivlist\vskip2pt}
\newenvironment{m-definition}{\vskip3pt\refstepcounter{msct}\trivlist \itemindent 0pt %
\item[\hskip\labelsep\bf Definition \themsct]\ignorespaces}{\endtrivlist\vskip2pt}
\newenvironment{m-notation}{\vskip3pt\refstepcounter{msct}\trivlist \itemindent 0pt %
\item[\hskip\labelsep\bf Notation \themsct]\ignorespaces}{\endtrivlist\vskip3pt}
\newenvironment{m-example}{\vskip3pt\refstepcounter{msct}\trivlist \itemindent 0pt %
\item[\hskip\labelsep\bf Example \themsct]\ignorespaces}{\endtrivlist\vskip3pt}
\newenvironment{m-remark}{\vskip3pt\refstepcounter{msct}\trivlist \itemindent 0pt %
\item[\hskip\labelsep\bf Remark \themsct]\ignorespaces}{\endtrivlist\vskip3pt}
\newenvironment{m-question}{\vskip3pt\refstepcounter{msct}\trivlist \itemindent 0pt %
\item[\hskip\labelsep\bf Question.]\ignorespaces}{\endtrivlist\vskip3pt}
\newenvironment{thm-nono}[1]{\vskip3pt\trivlist \itemindent 0pt %
\item[\hskip\labelsep\bf Theorem~#1]\it\ignorespaces}{\endtrivlist\vskip3pt}
\newenvironment{lm-nono}[1]{\vskip3pt\trivlist \itemindent 0pt %
\item[\hskip\labelsep\bf Lemma~#1]\it\ignorespaces}{\endtrivlist\vskip3pt}
\newenvironment{conj-nono}[1]{\vskip3pt\trivlist \itemindent 0pt %
\item[\hskip\labelsep\bf Conjecture~#1]\it\ignorespaces}{\endtrivlist\vskip3pt}
\newenvironment{m-thank}{\vskip3pt\trivlist \itemindent 0pt %
\item[\hskip\labelsep\it Acknowledgments]\ignorespaces}{\endtrivlist\vskip3pt}
\newenvironment{m-proof}{\vskip2pt\trivlist \itemindent 0pt %
\item[\hskip\labelsep\it Proof.]\ignorespaces}{\hfill$\Box$\endtrivlist\vskip3pt}%
\newenvironment{m-asmp}{\vskip3pt\trivlist \itemindent 0pt %
\item[\hskip\labelsep\bf Assumption.]\ignorespaces}{\hfill\endtrivlist\vskip3pt}%

\newcommand{\bibauth}[2]{\textrm{{#1}~{#2}}}
\newcommand{\bibtitl}[1]{\textit{#1}}
\newcommand{\bibjnyp}[4]{\textrm{#1} \textbf{#2} (#3), #4}
\newcommand{\bibinbook}[3]{In: \textrm{#1}\textrm{, #2}\textrm{, #3}}
\newcommand{\bibbook}[4]{\textit{#1}. {#2} {#3}, {#4}}

\let\lar\longrightarrow

\let\hra\hookrightarrow
\let\mt\mapsto

\font\tenmsa=msam10 %
\newcommand\hdashpiece{%
{\vrule height2.75pt depth-2.35pt width2.3pt \kern1.7pt}}%
\newcommand\hdashpieces{%
{\hdashpiece\hdashpiece\hdashpiece\hdashpiece}}%
\let\dashto\dashrightarrow
\newcommand\dashar{\mathrel{%
\hdashpieces\kern-0.4pt\hbox{\tenmsa K}}}%

\let\euf\EuScript 
\let\cal\mathcal
\let\mbb\mathbb

\let\mfrak\mathfrak
 
\DeclareFontFamily{OT1}{rsfs}{}
\DeclareFontShape{OT1}{rsfs}{n}{it}{<->rsfs10}{}
\DeclareMathAlphabet{\crl}{OT1}{rsfs}{n}{it}

\newcommand\uset[2]{{\disp\mathop{\mbox{$#2$}}_{#1}}}
\newcommand\oset[2]{{\disp\mathop{\mbox{$#2$}}^{#1}}}
\newcommand\ouset[3]{{\oset{#2}{\uset{#1}{#3}}}}

\let\unbar\underbar

\let\tld\tilde

\let\nit\noindent

\let\disp\displaystyle
\let\srel\stackrel
\let\vphi\varphi
\let\eps\epsilon

\newcommand\rd{{\rm d}} 

\newcommand\Aut{\operatorname{\textrm{Aut}\kern1pt}}
\newcommand\cAut{\operatorname{\mathcal{A}\kern-1pt\textit{ut}\kern1pt}}
\newcommand\End{\operatorname{\rm{End}\kern1pt}}
\newcommand\cEnd{\operatorname{\mathcal{E}\kern-1pt\textit{nd}\kern1pt}}

\newcommand\cHom{\operatorname{\mathcal{H}\kern-1pt\textit{om}\kern1pt}}

\newcommand\Ker{{\rm Ker}}

\newcommand\Pic{\mathop{\rm Pic}\nolimits}
\newcommand\Proj{\mathop{\rm Proj}\nolimits}

\newcommand\Spec{\mathop{\rm Spec}\nolimits}

\newcommand\Sym{\mathop{\rm Sym}\nolimits}
\newcommand\invq{{\slash\kern-2.5pt\slash}}
\newcommand\pr{\mathop{\rm pr}\nolimits}
\newcommand\rk{{\rm rk}}

\let\si\sigma
\let\Si\Sigma

\let\sm\setminus

\newcommand\bbk{\mbox{\rm I\kern-1.5pt k}}
\newcommand\sbbk{\hbox{\scriptsize I{\kern-.8pt}k}}

\newcommand\bone{{1\kern-0.57ex\rm l}}

\newcommand\eA{{\euf A}}

\newcommand\eE{{\euf E}}

\newcommand\eF{{\euf F}}
\newcommand\eG{{\euf G}}
\newcommand\eI{{\euf I}}

\newcommand\cL{{\cal L}}
\newcommand\eL{{\euf L}}

\newcommand{\eM}{{\euf M}}

\newcommand\eN{{\euf N}}
\newcommand\cO{{\cal O}}
\newcommand\eO{{\euf O}}

\newcommand\cY{{\cal Y}}

\newcommand\eT{{\euf T}}

\newcommand\codim{{\rm codim}}

\let\ges\geqslant
\let\les\leqslant

\newcommand\mx{{\rm max}}

\newcommand{\cst}[1]{\mathop{\rm ct}^{#1}\nolimits}

\newcommand{\Hilb}{\mathop{\rm Hilb}\nolimits}
\newcommand{\kk}{{\Bbbk}}
\newcommand{\lcit}{{\textit{loc.\,cit.}}}
\newcommand{\red}{{\rm red}}
\newcommand{\reg}{\mathop{\rm reg}\nolimits}
\newcommand{\trdeg}{\mathop{\rm trdeg}\nolimits}
\newcommand{\rst}{{\upharpoonright}}

\newcommand{\lbkt}{{[\kern-1.75pt[}}
\newcommand{\rbkt}{{]\kern-1.75pt]}}

\begin{document}

\title[Subvarieties with $q$-ample normal bundle]%
{Subvarieties with partially ample normal bundle}
\author{Mihai Halic}
\email{mihai.halic@gmail.com}
\address{CRM, UMI 3457, Montr\'eal~H3T~1J4, Canada}


\begin{abstract}
We show that local complete intersection subvarieties of smooth projective varieties, which have partially ample normal bundle, possess the G2-property. This generalizes results of Hartshorne and B\u{a}descu-Schneider. 
\end{abstract}

\maketitle


\section*{Introduction}

Hartshorne~\cite{hart-cdav,hart-as} investigated the cohomological properties of pairs $(X,Y)$, where $X$ is a projective scheme which is regular in a neighbourhood of a local complete intersection---lci, for short---subscheme $Y$ with ample normal bundle. He showed that, on one hand, $Y$ is G2 that is, the formal completion $\hat{X}_Y$ determines an \'etale neighbourhood of $Y$. On the other hand, the cohomology groups of coherent sheaves on the complement $X\sm Y$ are finite dimensional, above appropriate degrees. 

The ampleness of the normal bundle can be weakened. On the complex-analytic side, it suffices either a Hermitian metric with partially positive curvature (cf. \cite{griff,comm+graut}). On the algebraic side, B\u{a}descu-Schneider~\cite{bad+schn} addressed the \emph{globally generated}, partially ample case (in the sense of Sommese) by reducing the problem to~\cite{hart-cdav}. Their results mainly apply---due to the global generation of the normal bundle---to subvarieties of homogeneous varieties.  A comprehensive reference for the algebraic approach is B\u{a}descu~\cite{bad}. 

Subvarieties with $q$-ample normal bundle have not been investigated yet. Here we are referring to the cohomological partial ampleness~\cite{arap,tot}. It is less restrictive than Sommese's~\cite{soms} and also more flexible, being a numerical condition. There are numerous subvarieties with partially ample, but neither ample nor globally generated normal bundle. Their ubiquity is, in our opinion, a strong motivation to systematically study their properties. 

The main result of this article is stated below. It generalizes Hartshorne~\cite[Theorem 6.7]{hart-cdav}, B\u{a}descu-Schneider~\cite[Theorem 1]{bad+schn}, and strengthens as well the formality principle---for $Y$ lci rather than smooth---due to Griffiths, Commichau-Grauert, Chen~\cite{griff,comm+graut,chen}.

\begin{thm-nono}{\rm(cf.~\ref{cor:G2},~\ref{thm:chen})} 
Let $X$ be a smooth irreducible projective variety defined over an algebraically closed field of characteristic zero, and $Y$ a connected, lci subscheme, with $(\dim Y-1)$-ample normal bundle. Then $Y$ is G2 in $X$ and the formality principle holds for $(X,Y)$.
\end{thm-nono}

We conclude the article with applications. It is worth mentioning that Voisin's strongly movable subvarieties~\cite{voisin-coniv} have non-pseudo-effective co-normal bundle, hence they enjoy the G2-property (cf.~\ref{prop:1+}).


\section{Background material}\label{sct:def}

\begin{m-notation}\label{XYN} 
We work over an algebraically closed field $\kk$ of characteristic zero. Throughout the article, $\mfrak X$ is a connected, noetherian formal scheme, regular and projective over $\kk$; $X$ stands for an irreducible projective variety---that is, reduced and irreducible---over $\kk$. 

Let $Y$ be either a subscheme of definition of $\mfrak X$---it is projective---, or a closed subscheme of $X$; in the latter case, we suppose $X$ is non-singular along $Y$. Let $\dim Y$ be the maximal dimension of its components---we assume that all are at least $1$-dimensional---, $\codim_X(Y):=\dim X-\dim Y$ (if $Y\subset X$). Let $\eI_Y\subset\cal O_{\mfrak X}$ (resp. $\subset\eO_X$) be the sheaf of ideals defining $Y$; for $a\ges0$, $Y_{a}$ is the subscheme defined by $\eI_{Y}^{a+1}$. The formal completion of $X$ along $Y$ is $\hat X_Y:=\disp\varinjlim Y_a$; it is regular and projective. 

If $Y$ is lci in $\mfrak X$, we denote its normal sheaf by $\eN=\eN_Y:=(\eI_Y/\eI_Y^2)^\vee$; it is locally free of rank $\nu$. The structure sheaves of the various thickenings $Y_a$ fit into the exact sequences:
\begin{equation}\label{eq:Ya}
0\to\Sym^a(\eN^\vee)\to\eO_{Y_a}\to\eO_{Y_{a-1}}\to0,\;\forall a\ges1.
\end{equation}

For a coherent sheaf $\eG$, we denote $h^t(\eG):=\dim_{\kk}H^t(\eG)$; for a field extension $K\hra K'$, $\trdeg_KK'$ is the transcendence degree; $\cst{A,B,...}$ stands for a real constant depending on the quantities $A,B,\dots$. A \emph{line (resp. vector) bundle} is an \emph{invertible (resp. locally free) sheaf}.

We recall some terminology due to Hironaka-Matsumura~\cite{hir+mats}. Suppose $Y$ is connected; let $K(\hat X_Y)$ be the field of formal rational functions on $X$ along $Y$ (cf.~\cite[Lemma 1.4]{hir+mats}). 
\begin{itemize}[leftmargin=5ex]
\item[$\bullet$] $Y$ is G1 in $X$, if $H^0(\hat X_Y,\eO_{\hat X_Y})=\kk$;
\item[$\bullet$] $Y$ is G2 in $X$, if $K(X)\hra K(\hat X_Y)$ is finite;
\item[$\bullet$] $Y$ is G3 in $X$, if $K(X)\hra K(\hat X_Y)$ is an isomorphism.
\end{itemize}
\end{m-notation}


\subsection{Cohomological \textit{q}-ampleness}\label{ssct:cohom-q} 
This notion was introduced by Arapura and Totaro. 

\begin{m-definition}\label{def:q-line}
Let $Y$ be a projective scheme, $\eA\in\Pic(Y)$ an ample line bundle.
\begin{enumerate}[leftmargin=5ex]
\item[\rm(i)]  (cf. \cite[Theorem 7.1]{tot}) 
An invertible sheaf $\eL$ on $Y$ is \emph{$q$-ample} if, for any coherent sheaf $\eG$ on $X$, holds: 
$$\,\exists\,\cst{\eG}\;\forall\,a\ges \cst{\eG}\;\forall\,t>q,\; H^t(Y,\eG\otimes\eL^{a})=0.$$
\item[]
It's enough to verify the property for $\eG=\eA^{-k}, k\ges1$ (cf. \cite[Theorem 6.3, 7.1]{tot}).

\item[\rm(ii)]  (cf. \cite[Lemma 2.1, 2.3]{arap}) 
A locally free sheaf $\eE$ on $Y$ is \emph{$q$-ample} if $\eO_{\mbb P(\eE^\vee)}(1)$ on $\mbb P(\eE^\vee):=\Proj(\Sym^\bullet_{\eO_Y}\eE)$ is $q$-ample. It is equivalent saying that, for any coherent sheaf $\eG$ on $Y$, there is $\cst{\eG}>0$ such that: 
$$\,H^t(Y,\eG\otimes\Sym^a(\eE))=0,\;\forall t>q,\; \forall a\ges \cst{\eG}.$$ 
\item[] 
The \emph{$q$-amplitude of $\eE$}, denoted $q^\eE$, is the smallest integer $q$ with this property. 
Note that $\eE$ is $q$-ample if and only if so is $\eE_{Y_{\red}}$ (cf. \cite[Corollary 7.2]{tot}). Also, any locally free quotient $\eF$ of $\eE$ is still $q$-ample; indeed, $\eO_{\mbb P(\eF^\vee)}(1)=\eO_{\mbb P(\eE^\vee)}(1)\otimes\eO_{\mbb P(\eF^\vee)}$.

\item[\rm(iii)] 
For a coherent sheaf $\eG$ on $Y$, let $\reg^\eA(\eG)$ be its Castelnuovo-Mumford regularity with respect to $\eA$ and $\;\reg^\eA_+(\eG):=\mx\{1,\reg^\eA(\eG)\}.$
\end{enumerate}
\end{m-definition}

The $q$-amplitude enjoys \emph{uniformity} and \emph{sub-additivity} properties.

\begin{m-theorem}\label{thm:unif-q}\quad 
{\rm(i) (cf. \cite[Theorem 6.4, 7.1]{tot})}
Let $Y$ be a projective scheme, $\eA,\eL\in\Pic(Y)$. We assume that $\eA$ is sufficiently ample---Koszul-ample, cf. \cite[p.~733]{tot}---, and $\eL$ is $q$-ample. Then there are $\cst{\eA,\eL}_1,\cst{\eA,\eL}_2>0$, such that for any coherent sheaf $\eG$ on $Y$ holds:
$$
H^t(Y,\eG\otimes\eL^a)=0,\;\;\forall t>q,\,
\forall a\ges\cst{\eA,\eL}_1\cdot\reg^\eA_+(\eG)+\cst{\eA,\eL}_2.
$$

\nit{\rm(ii) (cf. \cite[Theorem 3.4]{tot})} 
If $H^0(\eO_Y)=\kk$ then, for a locally free sheaf $\eE$ and coherent sheaf $\eG$ on $Y$, one has 
$$\reg^\eA(\eE\otimes\eG)\les\reg^\eA(\eE)+\reg^\eA(\eG).$$ 
Hence it holds: $\reg^\eA_+(\eE\otimes\eG)\les\reg^\eA_+(\eE)+\reg^\eA_+(\eG).$
\end{m-theorem}

\begin{m-theorem}\label{thm:subadd-q}{\rm(cf. \cite[Theorem 3.1]{arap})}
Let $0\to\eE_1\to\eE\to\eE_2\to 0$ be an exact sequence of locally free sheaves on $Y$. Then it holds: $\;q^\eE\les q^{\eE_1}+q^{\eE_2}$. 
\end{m-theorem}

For products there is a better estimate. 

\begin{m-lemma}\label{lm:q-prod}
Let $X_1,X_2$ be irreducible projective varieties and $\eE_1,\eE_2$ locally free sheaves on them, respectively. Let $\eE_1\boxplus\eE_2$ be the direct sum of their pull-backs to $X_1\times X_2$. Then we have: 
$\;q^{\eE_1\boxplus\eE_2}\les\max\{q^{\eE_1}+\dim X_2,q^{\eE_2}+\dim X_1\}.$ 
\end{m-lemma}

\begin{m-proof}
Let $\eA_1,\eA_2$ be ample line bundles on $X_1,X_2$, respectively,  $\eA_1\boxtimes\eA_2$ the tensor product of their pull-backs. For $k\ges 1$, $t>\max\{q^{\eE_1}+\dim X_2,q^{\eE_2}+\dim X_1\},a\gg0,$ it holds:\\[1ex]
$\begin{array}{l}
H^t\big(X_1\times X_2,(\eA_1^{-k}\boxtimes\eA_2^{-k})\otimes\Sym^a(\eE_1\boxplus\eE_2)\big)
\\[1.5ex]
\null\kern1em=
\uset{\genfrac{}{}{0pt}{}{t_1+t_2=t,}{a_1+a_2=a}}{\bigoplus}
H^{t_1}\big(X_1,\eA_1^{-k}\otimes\Sym^{a_1}(\eE_1)\big)\otimes
H^{t_2}\big(X_2,\eA_2^{-k}\otimes\Sym^{a_2}(\eE_2)\big)=0. 
\end{array}$
\end{m-proof}

\begin{m-lemma}\label{lm:q-comp}
One has the equivalence:
\\[.5ex]\centerline{ 
$\eL\in\Pic(Y)$ is $q$-ample $\;\Leftrightarrow\;\eL\otimes\eO_{Y'}$ is $q$-ample, $\forall\,Y'\subset Y$ irreducible.
}
\end{m-lemma}

\begin{m-proof}
If $Y=Y'\cup Y''$ is the union of distinct closed subschemes, one has: 
$$
\begin{array}{l}
0\to\eO_Y\to\eO_{Y'}\oplus\eO_{Y''}\to\eO_{Y'\cap Y''}\to0,\\ 
0\to\eI_{Y'}\oplus\eI_{Y''}\to\eO_Y\to\eO_{Y'\cap Y''}\to0.
\end{array}
$$
Now tensor the exact sequences by $\eL^m\otimes\eO_Y(-k)$ and take their cohomology.
\end{m-proof}


\subsection{(dim\textit{Y}--1)-ample vector bundles on $Y$}\label{ssct:y-1} 

Subvarieties $Y\subset X$ with $(\dim Y-1)$-ample normal bundle will play an essential role. The following is analogous to Totaro's result for invertible sheaves. 

\begin{m-proposition}\label{prop:y-1}{\rm(cf. \cite[Theorem 9.1]{tot})}
Let $\eE$ be a locally free sheaf on an irreducible projective variety $Y$ (reduced, irreducible). The statements are equivalent: 
\begin{enumerate}[leftmargin=5ex]
\item[\rm(i)] 
$\eE$ is $(\dim Y-1)$-ample. 

\item[\rm(ii)] 
$\eO_{\mbb P(\eE)}(1)$ is not pseudo-effective, where $\mbb P(\eE):=\Proj(\Sym^\bullet\eE^\vee)$.
\item[] In this case, we say that \emph{$\eE^\vee$ is not pseudo-effective}. 

\item[\rm(iii)]  
There is a dominant morphism $\vphi:C_S\to Y$, with $S$ affine and $C_S$ an integral curve over $S$, such that the following conditions are satisfied: 
\begin{enumerate}
\item[\rm(1)] 
$\vphi^*\eE$ admits a line sub-bundle $\eM$ which is relatively ample for $C_S\to S$; 
\item[\rm(2)] 
Let $S_y\subset S$ be the curves passing through the general point $y\in Y$ and $\eM_{S_y}$ the restriction of $\eM$ to  $C_{S_y}$. 
\item[] 
Then the points $\{[\eM_{s,y}]\}_{s\in S_y}$, corresponding to $\eM_{s,y}$, cover an open subset of $\mbb P(\eE_y)$. (For shorthand, we say that $\eM\subset\vphi^*\eE$ is \emph{movable}.)
\end{enumerate}
\end{enumerate}
If $Y$ is reducible, the conditions {\rm(ii),~(iii)} must hold for all its irreducible components.
\end{m-proposition}

\begin{m-proof}
The last statement follows from~\ref{lm:q-comp}. Let $\eO_Y(1)$ be an ample line bundle on $Y$. Its dualizing sheaf $\omega_Y$ is torsion free of rank one, and $\eO_Y(-c)\subset\omega_Y\subset\eO_Y(c)$ for some $c>0$ (cf.~\cite[\S9]{tot}); hence the $(\dim Y-1)$-ampleness means:  
$H^0(Y,\omega_Y\otimes\eL\otimes\Sym^a\eE^\vee)=0,$ $\forall\eL\in\Pic(Y),$ $a>\cst{\eL}$. It is equivalent to $H^0(Y,\eM\otimes\Sym^a\eE^\vee)=0,$ $\forall\eM\in\Pic(Y),\, a>\cst{\eM}$, and to: 
$$
H^0(\mbb P(\eE),\eM\otimes\eO_{\mbb P(\eE)}(a))=0,\,\forall\eM\in\Pic(\mbb P(\eE)),\forall a>\cst{\eM}.
$$
The last condition is the $(\dim\mbb P(\eE)-1)$-ampleness of $\eO_{\mbb P(\eE)}(-1)$; (i)$\Leftrightarrow$(ii) follows. 

The equivalence (ii)$\Leftrightarrow$(iii) is the duality (cf.~\cite[Theorem 0.2]{bdpp}), for $X=\mbb P(\eE)$. However, \lcit\ requires $X$ to be smooth. Thus we must prove the following.

\nit\unbar{\textit{Claim}} Let $(X,\eO_X(1))$ be a $d$-dimensional projective variety, $\eL\in\Pic(X)$. It holds: 
$$
\eL\;\text{is $(d-1)$-ample}\;\Leftrightarrow\;
\exists\,\text{movable curve}\;C\to X\;\text{such that}\;\eL\cdot C>0.
$$
$(\Rightarrow)$ Let $\tld X\srel{\si}{\to}X$ be a desingularization of $X$ with exceptional locus $E$. Then $\si^*\eL$ is $(d-1)$-ample. Indeed, we may assume that $\tld\eA:=\big(\si^*\eO_X(1)\big)(-E)$ is ample on $\tld X$, hence:
$$
H^0(\tld X,\si^*\eL^{-m}\otimes\tld\eA^k)
\subset H^0(X,\eL^{-m}\otimes\eO_X(k)\otimes\si_*\eO_{\tld X})=0,\;k>0,\;m\gg\cst{k}.
$$
For the last step, $\si_*\eO_{\tld X}$ is torsion-free of rank one, so $\eO_X(-c)\subset\si_*\eO_{\tld X}\subset\eO_X(c)$ for an appropriate $c>0$. Then $\si^*\eL^{-1}$ is not pseudo-effective, so there is a movable curve $C\to\tld X$ such that $\si^*\eL\cdot C>0$.

\nit$(\Leftarrow)$  Let $C\to X$ be a movable curve and suppose $\eL$ is not $(d-1)$-ample. There is $k_0>0$ and a strictly increasing sequence ${\{m_t\}}_t\subset\mbb Z$, such that $H^0(X,\eL^{-m_t}\otimes \eA^{k_0})\neq0$. It follows: $0\les -m_t(\eL\cdot C)+k_0\eO_X(1)\cdot C,\;\forall\,t$, so $\eL\cdot C\les0$, a contradiction.
\end{m-proof}

Observe that the notion of pseudo-effective vector bundle used in~\cite[\S7]{bdpp} is more restrictive: it also requires that the projection of the non-nef locus of $\eO_{\mbb P(\eE)}(1)$ does not cover $Y$.



\section{Finite dimensionality results and the G2 property}\label{sct:fd+G2}

Hartshorne~\cite{hart-cdav} investigated the cohomological properties of lci subvarieties with ample normal bundle and of their complements. B\u{a}descu and Schneider \cite{bad+schn} extended his results to subvarieties with Sommese-$q$-ample (globally generated) normal bundle, hence their applications mainly concern homogeneous spaces. 

\subsection{Finite dimensionality}\label{ssct:finite-dim} 

The following generalizes results in~\cite[Section 5]{hart-cdav}.

\begin{m-theorem}\label{thm:finite-dim} 
Assume $Y$ is lci, let $q^\eN$ be the amplitude of its normal bundle $\eN$.
\begin{enumerate}[leftmargin=5ex]
\item[\rm(i)]
Consider $\cL\in\Pic(\mfrak{X})$, let $q^\cL$ be the amplitude of its restriction to $Y$. Let $\cal{F}$ be a locally free sheaf on $\mfrak{X}$, of finite rank. Then the following statements hold:
\begin{enumerate}
\item[\rm(a)] 
For $t<\dim Y-q^\eN$, $\,H^{t}(\mathfrak{X},\cal{F})$ is finite dimensional. In particular, if $q^\eN\leqslant\dim Y-1$ and $\mfrak{X}$ is connected, then $H^{0}(\mfrak{X},\cO_{\mfrak{X}})=\kk$.
\item[\rm(b)] 
$H^{t}(\mfrak{X},\cal{F}\otimes\cL^{-b})=0,
\text{ for }t<\dim Y-(q^\eN+q^\cL),\,b\gg0.$
\end{enumerate}
\item[\rm(ii)] 
Let $X$ be a projective scheme, non-singular along $Y$. Let $\eG$ be a coherent sheaf on $X\setminus Y$ and $\eL\in\Pic(X)$. The following statements hold: 
\begin{enumerate}
\item[\rm(a)]
$H^{t}(X\setminus Y,\eG)
\;\text{is finite dimensional,}\;t\ges\dim X-\dim Y+q^\eN,$
\item[\rm(b)] 
$H^{t}(X\setminus Y,\eG\otimes\eL^b)=0,\;t\ges\dim X-\dim Y+q^\eN+q^\eL,\;b\gg0.$
\end{enumerate}
\end{enumerate}
\end{m-theorem}

\begin{m-proof}
(i)(a) Use~\eqref{eq:Ya} and proceed as in \lcit,\ Theorem 5.1, Corollary 5.4. 

\nit(b)(cf. \lcit, Corollary 5.3) For $\eF:=\cal{F}\otimes\eO_{Y}$, $\eL:=\cL\otimes\eO_{Y}$, is enough to show:
$$\;
H^{t}(Y,\omega_Y\otimes\eF^\vee\otimes\Sym^{a}(\eN)\otimes\eL^b)=0,\quad\forall t>q^\eN+q^\cL,\forall a\ges0,b\gg0.
$$
But $\Sym^{a}(\eN)\otimes\eL^b$ is direct summand in $\Sym^{a+b}(\eN\oplus\eL)$, and $\eN\oplus\eL$ is $(q^\eN+q^\eL)$-ample. The vanishing holds for $a+b\ges\cst{\eF}$, \textit{e.g.} $a\ges 0$, $b\ges\cst{\eF}$.

\nit(ii) Use the formal duality~\cite[Theorem III.3.3]{hart-as} and the previous point.
\end{m-proof}


In  \cite[Expos\'e XIII, Conjecture 1.3]{groth}, Grothendieck discusses the finite dimensionality of the cohomology groups of coherent sheaves on the complement of lci subvarieties. Hartshorne addressed the issue for smooth subvarieties of projective spaces (cf. \cite[Corollary 5.7]{hart-cdav}). 

Let $S$ be a smooth projective variety and $E$ a principal $G$-bundle on it, with $G$ a connected linear algebraic group; let $P\subset G$ be a parabolic subgroup. Then $X:=E/P\srel{\pi}{\to}S$ is a locally trivial $G/P$-fibration. The co-ampleness ($ca$, for short) of homogeneous varieties has been explicitly computed by Goldstein~\cite{gold}. By definition, $q^{\eT_{G/P}}=\dim(G/P)-ca(G/P)$, hence $\eT_{X,\pi}:=\Ker(\rd\pi)$ is $q$-ample, for $q:=\dim X-ca(G/P)$. 

\begin{m-corollary}\label{cor:proj-bdl} 
Suppose $Y\subset X$ is a smooth $S$-family of subvarieties of relative codimension $\delta$, $\dim Y>\dim S$; that is, $\rd\pi_Y:\eT_Y\to\pi_Y^*\eT_S$ is surjective, $\codim_X(Y)=\delta$. Then $H^t(X\sm Y,\eG)$ is finite dimensional for $t\ges\delta+\dim X-ca(G/P)$, for all coherent sheaves $\eG$ on $X\sm Y$. 
\end{m-corollary}
Hartshorne's result corresponds to $S=\{\text{point}\}$, $G/P\cong\mbb P^n$, $t\ges\delta$.

\begin{m-proof} 
The exact diagram 
$$\xymatrix@R=1.35em{&0\ar[d]&0\ar[d]&&\\ 0\ar[r]&\eT_{Y,\pi_Y}\ar[r]\ar[d]&\eT_Y\ar[r]\ar[d]&\pi_Y^*\eT_S\ar[r]\ar@{=}[d]&0\\ 0\ar[r]&\eT_{X,\pi}\rst_Y\ar[r]\ar[d]&\eT_X\rst_Y\ar[d]\ar[r]&\pi^*\eT_S\rst_Y\ar[r]&0\\ &\eT_{X,\pi}\rst_Y\Big/\eT_{Y,\pi_Y}\ar@{=}[r]\ar[d]&\eN_{Y/X}\ar[d]&&\\ &0&0&&}$$
shows that $\eN_{Y/X}$ is a quotient of $\eT_{X,\pi}\rst_Y$, so is $q$-ample (cf.~\ref{def:q-line}(ii)); apply~\ref{thm:finite-dim}(ii). 
\end{m-proof}


\subsection{The G2 property}\label{ssct:G2} 
Here we generalize~\cite[Section 6]{hart-cdav}.  The difficulty to overcome is that several statements in there are proved for \emph{curves}, the general case being obtained by induction on the dimension.

\begin{m-lemma}\label{lm:L-ample}
{\rm(cf. \cite[Lemma 6.1]{hart-cdav})} 
Let $(Y,\eO_Y(1))$ be a projective scheme, $\eL\in\Pic(Y)$ and $\eE,\eF$ locally free sheaves on $Y$. Let 
$h_{\eF}(a,b):=h^{0}(Y,\eF\otimes \Sym^{a}(\eE^{\vee})\otimes\eL^{-b}),\,a,b\ges1.$ 
\begin{enumerate}[leftmargin=5ex]
\item[\rm(i)] 
If $\eL$ is $(\dim Y-1)$-ample, then it holds: 
\begin{equation}\label{eq:b>a}
h_{\eF}(a,b)=0,\;\text{for}\;\;
b\ges\cst{\eO_Y(1),\eL,\eE}_1\cdot\,a+\cst{\eO_Y(1),\eL,\eF}_2.
\end{equation}

\item[\rm(ii)] 
If $\eE$ is $(\dim Y-1)$-ample, then it holds: 
\begin{equation}\label{eq:a>b}
h_{\eF}(a,b)=0,\;\text{for}\;\;
a\ges\cst{\eO_Y(1),\eE,\eL}_1\cdot\,b+\cst{\eO_Y(1),\eE,\eF}_2.
\end{equation}
\end{enumerate}
\end{m-lemma}

\begin{m-proof} 
We fix $\eO_Y(1)$ sufficiently ample (cf.~\ref{thm:unif-q}) and consider the regularity with respect to it. Also, we may assume that $Y$ is irreducible; let $\omega_{Y}$ be its dualizing sheaf.\smallskip 

\nit(i) 
There is $c_{0}=c_0(Y)\geqslant1$ such that $\eO_Y(-c_0)\subset\omega_Y$, so it holds:
$$
\begin{array}{rl}
h^{0}(Y,\eF\otimes\Sym^{a}(\eE^{\vee})\otimes\eL^{-b})
&
\les h^{0}(Y,\omega_{Y}\otimes\eF(c_0)\otimes\Sym^{a}(\eE^{\vee})\otimes\eL^{-b})
\\[1.5ex]&
= 
h^{\dim Y}(Y,\eF^\vee(-c_0)\otimes\Sym^{a}(\eE)\otimes\eL^b).
\end{array}
$$
\nit\unbar{\textit{Claim}} The right hand side vanishes for $b$ as in \eqref{eq:b>a}. Indeed, we replace $\eF\rightsquigarrow\eF(-c_0)$ and verify the statement for $h^{\dim Y}(\eF^\vee\otimes \Sym^{a}(\eE)\otimes\eL^{b})$. The effect of the replacement is $\reg\eF^\vee\rightsquigarrow\reg\eF^\vee-c_0$, with $c_0$ depending on $Y$. 
Now observe that is enough to prove the claim for $Y$ reduced---so $H^0(\eO_Y)=\kk$---and for coherent sheaves $\eG$ on $Y$. 

Indeed, for  $\eI:=\Ker(\eO_Y\to\eO_{Y_{\text{red}}})$, there is $r>0$ such that $\eI^r=0$, so $\eO_Y$ admits a filtration (similar to \eqref{eq:Ya}) by the quotients $\eI^{k-1}/\eI^k$, $1\les k\les r$, which are $\eO_{Y_{\text{red}}}$-modules; now we may use the estimates for $\eF^\vee\otimes(\eI^{k-1}/\eI^k)$ on $Y_{\text{red}}$, which is coherent. Property~\ref{thm:unif-q} yields: 
$$
H^{\dim Y}(Y,\eG\otimes\Sym^a\eE\otimes\eL^b)=0,\;
\forall b\ges\cst{\eO_Y(1),\eL}_1\cdot\reg_+(\eG\otimes\Sym^a\eE)+\cst{\eO_Y(1),\eL}_2.
$$
But $\Sym^a\eE$ is a summand of $\eE^{\otimes a}$, so 
$\reg_+(\eG\otimes\Sym^a\eE)\les a\cdot\reg_+(\eE)+\reg_+(\eG),$ thus~\eqref{eq:b>a} holds for 
$b\ges\cst{\eO_Y(1),\eL}_1\cdot\,(a\cdot\reg_+(\eE)+\reg_+(\eG))+\cst{\eO_Y(1),\eL}_2.$\smallskip 

\nit(ii) 
We may assume that $Y$ is reduced. If $\eG$ is coherent on $Y$, 
$h^{\dim Y}(\eG\otimes \Sym^{a}(\eE)\otimes\eL^{b})$ vanishes for $a\ges\cst{\eO_Y(1),\eE}_1\cdot\reg_+(\eG\otimes\eL^b)+\cst{\eO_Y(1),\eE}_2$, and $\reg_+(\eG\otimes\eL^b)\les b\reg_+\eL+\reg_+\eG$. 
\end{m-proof}

\begin{m-proposition}\label{prop:L-any}
{\rm(cf. \cite[Theorem 6.2, Corollary 6.6]{hart-cdav})} 
Let the situation be as in \ref{XYN}. Suppose $Y$ is lci and its normal bundle $\eN$ is $(\dim Y-1)$-ample, of rank $\nu$. For any locally free sheaf $\cal F$ and invertible sheaf $\cal{L}$ on $\mfrak{X}$, there is a polynomial of degree $\dim Y+\nu$ such that: 
$$
h^{0}(\mfrak{X},\cal F\otimes\cal{L}^{b})\les P^{Y,\cal L,\cal F}_{\dim Y+\nu}(b),
\;\text{ for }b\gg0.
$$
\end{m-proposition}

\begin{m-proof} 
Let $\eA\in\Pic(Y)$ be sufficiently (Koszul) ample, such that $\eA^{-1}\subset\omega_Y$; denote $\eF:=\cal F\otimes\eO_Y,\eL:=\cal L\otimes\eO_Y$. For $\gamma:=\cst{\eA,\eN,\eL}_1+1, b>\cst{\eA,\eN,\eF}_2$ (cf.~\eqref{eq:a>b}), it holds:
$$
h^{0}(\mfrak{X},\cal F\otimes\cal{L}^{b})
\leqslant
\ouset{a=0}{\infty}{\sum} h^{0}(Y,\eF\otimes\Sym^{a}(\eN^{\vee})\otimes\eL^{b})
=
\ouset{a=0}{\gamma b}{\sum} h^{0}(Y,\eF\otimes\Sym^{a}(\eN^{\vee})\otimes\eL^{b}).
$$
Since $\eF\subset(\eA^{c_0})^{\oplus\rk\eF}, c_0=\reg^\eA_+\eF^\vee$, it is enough to consider $\eF=\eA^{c_0}$. 

Consider $S:=\mathbb{P}(\eO_{\mathbb{P}(\eN^{\vee})}(-1)\oplus\eO_{\mathbb{P}(\eN^{\vee})})$ and $\eO_{S}(1)$ the relatively ample invertible sheaf on it. The right hand side above can be re-written:
\begin{longtable}{l}
rhs=
$
\ouset{a=0}{\gamma b}{\sum} h^{0}(Y,\eA^{c_0}\otimes\Sym^{a}(\eN^{\vee})\otimes\eL^{b})
\les
\ouset{a=0}{\gamma b}{\sum} 
h^0(Y,\omega_{Y}\otimes\eA^{c_0+1}\otimes\Sym^{a}(\eN^\vee)\otimes\eL^{b})$
\\[2ex]
$=
\ouset{a=0}{\gamma b}{\sum} 
h^{\dim Y}(Y,\eA^{-c_0-1}\otimes\Sym^{a}(\eN)\otimes\eL^{-b})$
\\[2ex] 
$=
\ouset{a=0}{\gamma b}{\sum}
h^{\dim Y}(\mathbb{P}(\eN^{\vee}),\eA^{-c_0-1}\otimes\eO_{\mathbb{P}(\eN^{\vee})}(a) \otimes\eL^{-b})
=
h^{\dim Y}(S,\eA^{-c_0-1}\otimes\eO_{S}(\gamma b)\otimes\eL^{-b}).$
\end{longtable}
\nit But $h^{\dim Y}(S,\eO_{S}(\gamma b)\otimes\eL^{-b})$ is dominated by a polynomial in $b$, depending on $\eO_S(\gamma)\otimes\eL^{-1}$, of degree at most $\dim S=\dim Y+\nu$ (cf. \cite[1.2.33]{laz}). To include $\eA^{-c_0-1}$, use 
$$
0\to\eA^{-c_0-1}\to\eO_Y\to\eO_{Y_1}\to0,\;\;\dim Y_1=\dim Y-1,
$$
which yields: 
$\text{rhs}\les h^{\dim X}(\eO_S(\gamma b)\otimes\eL^{-b}) + h^{\dim X-1}(\eO_S(\gamma b)\otimes\eL^{-b}\rst_{Y_1}).$ 
\end{m-proof}
With these preparations, the proof of the following theorem is identical to \lcit 

\begin{m-theorem}\label{thm:G2}
{\rm(cf. \cite[Theorem 6.7]{hart-cdav})} 
Let the situation be as in \ref{XYN}. We assume: 
\begin{enumerate}[leftmargin=5ex]
\item[$\bullet$] $Y$ is connected, lci, $\dim Y\geqslant1$;
\item[$\bullet$] the normal bundle $\eN$ of $Y$ is $(\dim Y-1)$-ample.
\end{enumerate}
Then the following statements hold: 
\begin{enumerate}[leftmargin=5ex]
\item[\rm(i)] 
$\trdeg_{\kk}K(\mfrak{X})\leqslant\dim Y+\rk\eN$;
\item[\rm(ii)] 
If $\trdeg_{\kk}K(\mfrak{X})=\dim Y+\rk\eN$, then $K(\mfrak{X})$ is a finitely generated extension of $\kk$.
\end{enumerate}
\end{m-theorem}

\begin{m-corollary}\label{cor:G2}
{\rm(cf. \cite[Corollary 6.8]{hart-cdav})} 
Let $X$ be a projective scheme, non-singular in a neighbourhood of a closed, connected, lci subscheme $Y$ with $(\dim Y-1)$-ample normal bundle. Then $Y$ is G2 in $X$. 
\end{m-corollary}

\begin{m-proof}
Indeed, $K(X)$ is a subfield of $K(\hat X_Y)$, so $\trdeg_\kk K(\hat X_Y)\ges\dim X=\dim Y+\nu$. Hence we are in the case (ii) of the previous theorem. 
\end{m-proof}
The result is optimal, one can not conclude that $Y$ is G3 (cf.~\cite[Example p.~199]{hart-as}.


\subsection{A formality criterion}\label{ssct:chen}

One says that the \emph{formal principle} holds for a pair $(X,Y)$ consisting of a scheme $X$ and a closed subscheme $Y$ if the following condition is satisfied: for any other pair $(Z,Y)$ such that $\hat Z_Y\cong\hat X_Y$, extending the identity of $Y$, there is an isomorphism between \'etale neighbourhoods of $Y$ in $X$ and in $Z$ which induces the identity on $Y$.

\begin{m-theorem}\label{thm:chen} 
In the situation~\ref{cor:G2}, the formal principle holds for $(X,Y)$. 
\end{m-theorem}
This simplifies and strengthens~\cite[Theorem 3]{chen}, since $Y$ is only lci, rather than smooth. 

\begin{m-proof}
Corollary \ref{cor:G2} implies that $Y$ is G2 in $X$. But, in this case, Gieseker proved (cf. \cite[Theorem 4.2]{gskr}, \cite[Corollary 9.20, 10.6]{bad}) that the formality holds for $(X,Y)$. 
\end{m-proof}

There are similar results in complex analytic setting. Griffiths~\cite{griff} investigated the formality/rigidity of smooth subvarieties $Y\subset X$ whose normal bundle $\eN_{Y/X}$ admits a Hermitian metric with curvature of signature $(s,t)$, $s+t=\dim Y$, and proves in [\textit{ibid.}, II.~\S2,\,3] the rigidity of the embedding for $s\ges 2$. The main cohomological property of vector bundles admitting metrics of curvature with mixed signature $(s,t)$ is that of being $(\dim Y-s)$-ample (cf. \cite[Proposition 28, p.~257]{andr+grau}, \cite[(7.28), p.~432]{griff}). 

On the other hand, Commichau-Grauert \cite[Satz 4]{comm+graut} proved the formality for subvarieties with $1$-positive normal bundle. Note that a $1$-positive vector bundle on a smooth projective variety $Y$ is $(\dim Y-1)$-ample (cf. \cite[Satz 2]{comm+graut}). 

We conclude that the cohomological approach adopted in this article yields under weaker assumptions the rigidity results obtained in \cite{griff,comm+graut}.


\section{Examples of subvarieties with partially ample normal bundle}
\label{sct:apply}

In this section we assume that $X$ is a smooth projective variety. 

\subsection{Elementary operations}\label{ssct:elem}

\begin{m-corollary}\label{cor:elem} 
\begin{enumerate}[leftmargin=5ex]
\item[\rm(i)] 
Let $Y_2\subset Y_1\subset X$ be connected lci, $\dim Y_2\ges1$. Suppose $\eN_{Y_2/Y_1}, \eN_{Y_1/X}$ are respectively $q_2$-, $q_1$-ample, with $q_1+q_2<\dim Y_2$. Then $Y_2$ is G2 in $X$.
\smallskip 

\item[\rm(ii)] 
Suppose $Y_1, Y_2$ are lci in $X$, 
$
\codim (Y_1\cap Y_2)=\codim (Y_1)+\codim(Y_2),
$
and $\eN_{Y_j/X}$ is $q_j$-ample, for $j=1,2$. Then $\eN_{Y_1\cap Y_2/X}$ is $(q_1+q_2)$-ample. 
\smallskip

\item[\rm(iii)] 
Suppose $Y_j\subset X_j$ are connected lci and $\eN_{Y_j/X_j}^\vee$ is not pseudo-effective (so $Y_j\subset X_j$ is G2), for $j=1,2$. Then $Y_1\times Y_2$ is lci and G2 in $X_1\times X_2$. 
\smallskip

\item[\rm(iv)] 
Let $f:X'\to X$ be a surjective, flat morphism. Suppose $Y\subset X$ is lci and $\eN_{Y/X}$ is $(\dim Y-1)$-ample. Then $Y':=f^{-1}(Y)\subset X$ is lci and $\eN_{Y'/X'}$ is $(\dim Y'-1)$-ample.
\end{enumerate}
\end{m-corollary}

\begin{m-proof}
(i)-(iii) are consequences of the sub-additivity~\ref{thm:subadd-q} and \ref{lm:q-prod}, applied to appropriate normal bundle sequences. 
For (iv), note that $f$ is equidimensional, $\eN_{Y'/X'}=f^*\eN_{Y/X}$. Apply Leray's spectral sequence to $Y'\to Y$.
\end{m-proof}

\begin{m-corollary}\label{cor:G2-inters}
Suppose $Y_1, Y_2$ are lci in $X$, 
$\codim (Y_1\cap Y_2)=\codim (Y_1)+\codim(Y_2),$ and $\eN_{Y_j/X}$ is $q_j$-ample, for $j=1,2$. If $Y_1\cap Y_2$ is connected and $q_2<\dim(Y_1\cap Y_2)$, e.g. $q_2=0$, then $Y_1\cap Y_2$ is G2 in $Y_1$. 
\end{m-corollary}

\begin{m-proof}
Note that $\eN_{Y_1\cap Y_2/Y_1}\cong{\eN_{Y_2/X}\rst}_{Y_1\cap Y_2}$. 
\end{m-proof}


\subsection{Strongly movable subvarieties}\label{ssct:strg-mov}
(cf.~\cite[Section~2]{voisin-coniv})~A class of examples of subvarieties having the G2-property are the strongly movable subvarieties introduced by Voisin~\cite[Section~2]{voisin-coniv}, in the attempt to geometrically characterize big subvarieties. 

\begin{m-notation}\label{not:syx}
Let $\cY\srel{(\pi,\rho)}{\subset} S\times X$ be a flat family of lci subschemes of $X$, with $\rho$ dominant; then $\rho(\cY)$ contains an open subset $O$ of $X$. We may (and do) assume that $S,\cY$ are reduced, since so is $X$. The incidence variety $\Si$ is the component of $(\pi_1,\pi_2)(\cY\times_X\cY)\subset S\times S$ containing the diagonal; $\pi$ is a proper, so $\Si$ is closed. One obtains the Cartesian diagram:
\begin{equation}\label{eq:rho}
\xymatrix@R=1.35em@C=5em{
\cY_{\Si_o}\ar[r]\ar[d]^-{\pi}&
\cY_\Si\ar[r]\ar@/^2.5ex/[rr]^-{\rho_\Si}\ar[d]&
\cY\times\cY\ar[r]_-{(\pi_1,\rho_2)}\ar[d]_-{(\pi_1,\pi_2)}&S\times X.
\\ 
\Si_o\ar[r]&
\Si\ar[r]^-{\iota}&S\times S&
}
\end{equation}
For $o\in S$, denote $\Si_o:=\iota^{-1}(\{o\}\times S)$ and $\rho_o:=\rho_\Si\rst_{\Si_o}$. 
\end{m-notation}

\begin{m-definition}\label{def:mov+conn}
Suppose the general member of $\cY$ is irreducible. We say that the family $\cY$ is \emph{strongly movable}, if $\rho_\Si$ is dominant; then $\cY_{\Si_o}\srel{\rho_o}{\to}X$ is dominant, for $o\in S$ general, and $Y_o$ is strongly movable. An \emph{arbitrary} family $\cY$ is \emph{strongly movable} if so is its general member $Y_o$; that is, \emph{all} the irreducible components of $Y_o$ are strongly movable.
\end{m-definition}

\begin{m-proposition}\label{prop:1+}
Let $\cY$ be as above, $o\in S$ a non-singular point such that $Y_o$ is strongly movable. Then $\eN_{Y_o/X}$ is $(\dim Y_o-1)$-ample. Hence, if $Y_o$ is connected, it is G2 in $X$. 
\end{m-proposition}

\begin{m-proof} 
Let $\eT_{S,o}\srel{\delta}{\to} H^0(Y_o,\eN_{Y_o/X})$ be the infinitesimal deformation homomorphism. By~\ref{prop:y-1}, it is enough to prove that the restriction of $\eN_{Y_o/X}$ to the irreducible components of $Y_o$ are $(\dim Y_o-1)$-ample. Recall that $\eN_{Y_o/X}$ is $(\dim Y_o-1)$-ample if and only if so is its restriction to $Y_{o,\red}$. For $\xi\in\eT_{S,o}$, we denote $\hat v_\xi\in H^0(Y_{o,\red},\eN_{Y_o/X}\rst_{Y_{o,\red}})$ the restriction of $\delta(\xi)$ to $Y_{o,\red}$. Henceforth, we replace $Y_o$ by an irreducible component.  

We must find a movable morphism $C\srel{\vphi}{\to}Y_{o,\red}$, an ample line bundle $\eL_C\in\Pic(C)$, and a movable homomorphism $\eL_C\to\vphi^*\eN_{Y_o/X}$. We restrict ourselves to $\xi\in\eT_{\Si_{o},o}\subset\eT_{S,o}$. 

\nit\underbar{\textit{Claim~1}}\quad 
The vanishing locus of $\hat v_\xi$ is a non-empty, proper subset of $Y_{o,\red}$; for $\xi\in\eT_{\Si_o,o}$ variable, the vanishing loci of $\hat v_\xi$ cover an open subset of $Y_{o,\red}$. \\ 
The vector $\xi$ is determined by an arc $\Spec\big(\kk\lbkt\eps\rbkt\big)\srel{h}{\to}\Si_o\,$ through $o$. The defining property of $\Si_o$ implies that $h(\eps)=y_\eps\in Y_o\cap Y_{h(\eps)}$, $h(0)=y\in Y_o$. Since $Y_o$ is deformed at $y$ in a tangential direction, we deduce $\delta(\xi)_y=0$, so $\hat v_{\xi,y}=0.$

We claim that $\hat v_\xi\neq0$, for generic $\xi$, and their vanishing loci contain an open subset of $Y_{o,\red}$. Indeed, $\cY_{\Si_o}\srel{\rho_o}{\to}X$ is dominant, so $\eT_{\Si_o,o}\srel{\rd\rho_{o,y}}{\lar}\eN_{Y_{o,\red}/X,y}$ is surjective at a generic (smooth) point $y\in Y_{o,\red}$. Since $Y_o$ is lci, a computation in local coordinates shows that there is a non-trivial homomorphism $\eN_{Y_o/X,y}\to\Sym^k(\eN_{Y_{o,\red}/X,y})$, for some $k>0$ (\textit{e.g.} $k=1$, if $Y_o$ is reduced at $y$), hence $\hat v_{\xi,y}\neq0$. Second, $Y_o$ is strongly movable, so the points $y\in Y_o$ where $\delta(\xi)_y=0$, for some $\xi\in\eT_{\Si_o,o}$, cover an open subset of $Y_o$. 

\nit\underbar{\textit{Claim~2}}\quad 
Let $C\subset Y_o$ be a complete intersection curve which intersects the zero locus of $\hat v_\xi$ properly. By Claim~1, such curves are movable. Moreover, $\hat v_\xi$ extends to a pointwise injective homomorphism $\eL_C\subset\eN_{Y_o/X}\rst_C$, where $\eL_C$ is an ample line bundle.  The latter is movable too, because $\rd\rho_{o,y}$ is surjective at the generic point $y\in Y_o$. This is formalized as follows. 

Let $\cal C_R\srel{f}{\to}X$ be a movable curve ($R$ is a parameter variety). Consider the diagram 
$$
\xymatrix@R=1.5em@C=3em{
g^*\pr^*\eN\ar[d]&\pr^*\eN\ar[d]
&\eN:=\eN_{Y_o/X}\rst_{Y_{o,\red}}\ar[d]
\\
\eT_{\Si_o,o}\times\cal C_R\ar[r]^-{g=(\delta,f)}\ar@/^3ex/[u]^{g^*v}
&H^0(Y_{o,\red},\eN)\times Y_{o,\red}\ar[r]^-\pr\ar@/^3ex/[u]^{v}
&Y_{o,\red}
}
$$
where $v$ is the evaluation map. We may suppose that  $\cal C_R$ is such that its generic member intersects non-trivially and properly the zero locus of $v$. Note that $g^*v$ yields a rational map $\eT_{\Si_o,o}\times\cal C_R\dashto g^*\pr^*\mbb P(\eN)$ which extends to a morphism outside a closed subscheme $Z$ of codimension at least two. Its projection does not cover $\eT_{\Si_o,o}\times R$, hence we obtain a movable, relatively ample $\eL\subset g^*\pr^*\eN$. 
\end{m-proof}


\subsection{Varieties whose cotangent bundle is not pseudo-effective}
\label{ssct:not-pseff} 

Let $Y\subset X$ be a smooth subvariety, so $\eN_{Y/X}$ is a quotient of $\eT_X\rst_Y$.

\begin{m-notation}
For shorthand, denote $\mbb P:=\mbb P(\eT_X)$ and $\mbb P_Y:=\mbb P(\eT_X\rst_Y)$ its restriction to $Y$; let $\pi:\mbb P\to X$ be the projection. 
Define ${\rm Mov}(\mbb P_Y)_{\mbb Q}\subset H_2(\mbb P_Y;\mbb Q)$ to be the cone generated by the classes of movable curves on $\mbb P_Y$ and ${\rm Mov}(\mbb P)_{\mbb Q}$ similarly.
\end{m-notation}

\begin{m-corollary}\label{cor:mov}
Let $Y\subset X$ be a smooth subvariety such that $\eO_{\mbb P_Y}(1)$ is not pseudo-effective. Then $\eN_{Y/X}$ is $(dim Y-1)$-ample, so $Y$ is G2.
\end{m-corollary}

\begin{m-proof}
Theorem~\ref{thm:G2} applies, since $\eN_{Y/X}$ is $(\dim Y-1)$-ample.
\end{m-proof}

By using~\ref{prop:y-1} we are going to show that, for $Y\subset X$ sufficiently general, the partial ampleness of $\eT_X\rst_Y$ implies the non-pseudo-effectiveness of the cotangent bundle of $X$. The latter is a numerical condition/restriction on the ambient variety. Examples include rationally connected varieties---see below---and, possibly, Calabi-Yau varieties (cf.~\cite[Corollary 6.12]{dps}). 

\begin{m-lemma}\label{lm:mov} 
Let $Y\srel{\iota}{\to} X$ be a subvariety, $\dim Y>0$, and $H_2(\mbb P_Y;\mbb Q)\srel{\iota_*}{\to} H_2(\mbb P;\mbb Q)$ be the induced homomorphism. In the situations enumerated below, it holds: 
\begin{equation}\label{eq:mov}
\iota_*\big({\rm Mov}(\mbb P_Y)\big)_{\mbb Q}\subseteq{\rm Mov}(\mbb P)_{\mbb Q}.
\end{equation}
\begin{enumerate}[leftmargin=5ex]
\item[\rm(i)] 
An algebraic group $G$ acts on $X$ with an open orbit $O$, such that the stabilizer of a point $x\in O$ acts with open orbit on $\eT_{X,x}$, and $Y\cap O\neq\emptyset$. 
\item[\rm(ii)] 
$\kk$ is uncountable, and $Y$ is a  very general member of a dominant family. 
\end{enumerate}
Hence, if $\eO_{\mbb P_Y}(1)$ is not pseudo-effective, then $\eO_{\mbb P}(1)$ is the same.
\end{m-lemma}

\begin{m-proof} 
(i) The $G$-translates of a movable curve on $\mbb P_Y$ cover an open subset of $\mbb P(\eT_X)$. 

\nit(ii) Let $Y_S\subset S\times X$ be an $S$-flat family of subvarieties, with $S$ affine, dominating $X$. Curves on $X$ are parametrized by their Hilbert polynomial, of degree one, with integer coefficients. Let $\Pi$ be the countable set of polynomials occurring for movable curves on $Y_s,s\in S$. 

For $P\in\Pi$, denote $\Hilb^P_{Y_S/S}\srel{\pi_P\,}{\to}S$ the corresponding scheme. We are interested in the components corresponding to curves. For $s\in S$, let $\Pi_s\subset \Pi$ be the set of polynomials $P_s$ such that $\pi_{P_s}$ is not dominant; $\Pi_{\text{rigid}}:=\uset{s\in S}{\bigcup}\Pi_s$. The image of $\Hilb^{\Pi_{\text{rigid}}}_{Y_S/S}\to S$ is a countable union of proper subvarieties. 
Take $s'\in S$ in the complement ($\kk$ is uncountable); let $P_{s'}$ be the Hilbert polynomial of some movable curve $C_{s'}\subset Y_{s'}$. Then $P_{s'}\not\in\Pi_{\text{rigid}}$, so $\Hilb^{P_{s'}}_{Y_S/S}\srel{\pi_{P_{s'}}\;}{\lar}S$ is surjective. Let $\Pi':=\Pi\sm\Pi_{\text{rigid}}$. The components of $\Hilb^{\Pi'}_{Y_S/S}$ (corresponding to movable curves) dominate $S$, so they are flat over the very general point $o\in S$. 

We claim that movable curves on $Y_o$ are movable on $X$. Indeed, for $P_o$ as above, consider the universal curve $\cal C_S\subset\Hilb^{P_o}_{Y_S/S}\times_S Y_S$. The family $\cal C_{o}\subset \Hilb^{P_o}_{Y_o}\times Y_o$ dominates $Y_o$. By the continuity of $\Hilb^{P_o}_{Y_S/S}\times_S Y_S\to Y_S$, the same holds for $\cal C_{s}\subset \Hilb^{P_o}_{Y_s}\times Y_s$, for $s$ near $o\in S$. Finally, $Y_S\to X$ is dominant, so $\cal C_S$ covers an open subset of $X$. 
\end{m-proof}

\begin{m-lemma}\label{lm:rtl-conn}
Let $X$ be a smooth rationally connected variety. Then $\eO_{\mbb P(\eT_X)}(1)$ is not pseudo-effective, so $\eT_X$ is $(\dim X-1)$-ample.
\end{m-lemma}

\begin{m-proof}
Consider a very free rational curve: a dominant morphism $\mbb P^1\times S\srel{\vphi}{\to} X$, where $S$ is a variety, such that 
$\vphi^*\eT_X\rst_{\mbb P^1\times\{s\}}\cong\eO_{\mbb P^1}(1)\otimes\eG,\,\forall s\in S,$ with $\eG$ globally generated. A nowhere vanishing section $g_s\in H^0(\mbb P^1\times\{s\},\eG)$ yields the inclusion $j_{g_s}:\eO_{\mbb P^1}(1)\to\vphi^*\eT_X\rst_{\mbb P^1\times\{s\}}$; we still denote by $j_{g_s}$ the morphism $\mbb P^1\to\mbb P(\eT_X)$. Since $\eG$ is globally generated, there is $\tld S\subset S\times H^0(\eG)$ open, such that 
$\tld S\to{\rm Morphisms}(\mbb P^1,\mbb P(\eT_X)),\;\{s\}\times\{g_s\}\mt j_{g_s},$ yields a movable rational curve $\mbb P^1\times\tld S\srel{\tld\vphi}{\to}\mbb P(\eT_X)$ satisfying 
$\tld\vphi^*\eO_{\mbb P(\eT_X)}(-1)\rst_{\mbb P^1\times\{\tld s\}}\cong \eO_{\mbb P^1}(1),\;\forall \tld s\in\tld S$.
\end{m-proof}

It is well-known---particularly for projective spaces---that the G3 property of the diagonal  $\Delta_X:=\{(x,x)\mid x\in X\}\subset X\times X$ imply important connectedness results for the intersections of subvarieties in $X$ (cf.~\cite[Ch.~11]{bad}). Our results yield the G2-property of the diagonal; obviously, it is less than the G3-property, but it holds for a larger class of varieties. 

\begin{m-proposition}\label{prop:diag-G2}
Let $X$ be a smooth projective variety, whose cotangent bundle is not pseudo-effective (e.g. rationally connected). Then the diagonal $\Delta_X$ is G2 in $X\times X$. 
\end{m-proposition}

\begin{m-proof}
The normal bundle of $\Delta_X$ is isomorphic to $\eT_X$; we conclude by~\ref{cor:G2}.
\end{m-proof}


\end{document}